\documentclass[12pt,a4paper]{amsart}
\usepackage{amssymb,amsmath,latexsym}
\usepackage{epsfig}

\newtheorem{thm}{Theorem}[section]
\newtheorem{dfn}[thm]{Definition}
\newtheorem{cor}[thm]{Corollary}

\newtheorem{lemma}[thm]{Lemma}

\newcommand{\mc}[1]{\mathcal{#1}}
\newcommand{\nat}{\mathbb{N}}
\newcommand{\real}{\mathbb{R}}

\newfont{\menutt}{cmtt8}

\title[Flag Matroids]{A natural family of flag matroids}
\date{\today}
\author{Anna de Mier}
\address{Mathematical Institute\\ University of Oxford\\ 24--29 St Giles\\
 Oxford OX1 3LB\\
 United Kingdom} \email{ademier@gmail.com}

\begin{document}

\begin{abstract}
A flag matroid can be viewed as a chain of matroids linked by
quotients. Flag matroids, of which relatively few interesting
families have previously been known,  are a particular class of
Coxeter matroids. In this paper we give a family of flag matroids
arising from an enumeration problem that is a generalization of
the tennis ball problem. These flag matroids can also be defined
in terms of lattice paths and they  provide a generalization of
the lattice path matroids of [Bonin et al., JCTA 104 (2003)].
\end{abstract}

\maketitle

\section{Introduction and preliminaries}\label{sec:intro}

Flag matroids are a subclass of Coxeter matroids, but they can also be
described in pure matroid-theoretical terms. Roughly speaking, a flag matroid
is a collection of matroids on the same ground set that form a chain in the
strong order (i.e., they are quotients of each other). Flag matroids play
an important role in the theory of Coxeter matroids and also shed light on
ordinary matroid theory.
Nevertheless, not many classes of flag matroids have been studied up to now.
The goal of this paper is to introduce a new family of flag matroids based on
an enumeration problem and show how these flag matroids can be interpreted in
terms of lattice paths. We refer to~\cite{coxeter}, especially
to Chapter~1, for an introduction to flag matroids and the ideas behind them.

We assume the reader is familiar with the basic concepts of matroid
theory; we
follow the notation of Oxley's book~\cite{oxley}. We recall here only the
notion of quotient. Given two matroids $M$ and $N$ on the same ground set, $M$ is
a \emph{quotient} of $N$ if every flat of $M$ is a flat of $N$ (one can also say
that $M$ is a strong map image of $N$). In this case,
the rank of $M$ is at most the rank of $N$, with equality holding if and only
if $M$ and $N$ are equal.

We also need to say a few words about lattice path matroids. We
do not need lattice path matroids in general as defined in~\cite{lpm1},
but
only the subclass of nested matroids. These matroids have independently
arisen several times in the literature since at least 1965, and have been given a variety of names;
see~\cite{lpm2,lpm1} and the references therein
for definitions and results (in these papers, nested matroids are called
``generalized Catalan matroids'').

Let $P$ be a lattice path from $(0,0)$ to $(m,r)$ with steps $E=(1,0)$ and
$N=(0,1)$. Let $\mathcal{P}$ be the set of paths from $(0,0)$ to $(m,r)$ with
steps $E$ and $N$ and that do not go above $P$.
 For each path $Q \in \mathcal{P}$,
let $Q_N=\{i: \mbox{ step } i \mbox{ in } Q \mbox{ is } N\}$. We denote
by $[n]$ the set $\{1,2,\ldots, n\}$.

\begin{thm}\label{thm:defnested}
 The set  $\{Q_N: Q\in \mathcal{P}\}$ is the collection of bases of a matroid
$M[P]$ on the ground set $[m+r]$.
\end{thm}

A matroid
is \emph{nested} if it is isomorphic to $M[P]$ for some path $P$.
Hence, the bases of a nested matroid are in bijection with the lattice
paths that do not go above a certain fixed path $P$ (see the left side
of Figure~\ref{fig:tbpmat} for an example of a nested matroid on the set
 $[15]$; the path
highlighted corresponds to the basis $\{2,5,8,11,12,15 \}$).
The name nested comes from the fact that a nested matroid can also be defined
as a transversal matroid whose presentation consists of nested sets, and
also because of the following characterization of nested matroids in terms of cyclic
flats (recall that a flat is \emph{cyclic} if it is a union of circuits).

\begin{thm}\label{thm:nested}
A matroid is nested if and only if its cyclic flats form a chain under
 inclusion. Furthermore, the proper non-trivial cyclic flats of the matroid
 $M[P]$ are the
initial segments $[t]$ of $[m+r]$, where $t$ is such that step
$t$ of $P$ is $E$ and step $t+1$ is $N$.
\end{thm}

Our view on flag matroids is slightly different from that of~\cite{coxeter},
but it is easy to see that the two perspectives are equivalent.
The definition in~\cite{coxeter} is in terms of flags of sets, whereas ours relies on what
we call ordered partitions of a set. For the reader already familiar with
the theory of flag matroids, changing from one definition to the other should
be straightforward.

\begin{dfn}\label{dfn:orderedpartition}
An \emph{ordered $k$-partition} of a set $S$ is a $k$-tuple $(A_1,\ldots,A_k)$ of
non-empty sets with
$A_1\cup A_2 \cup \cdots \cup A_k=S$  and $A_i\cap A_j=
\emptyset$ whenever $i\neq j$. For  positive integers
$r_1,\ldots, r_k$ such that $r_1+\cdots +r_k=|S|$,
an \emph{$(r_1,r_2,\ldots, r_k)$-partition}
of $S$ is an ordered $k$-partition $(A_1,\ldots,A_k)$ of $S$ such that
$|A_i|=r_i$ for all $i$ with $1\leq i \leq k$.
\end{dfn}

The bases of a matroid $M$ on a set $S$ trivially determine a collection of
ordered  $2$-partitions of $S$: take all pairs of the form $(B, S-B)$, where
$B$ is a basis of $M$. The first axiom for flag matroids generalizes this idea; the
other two axioms arise from the definition of a flag matroid in terms of Coxeter
groups (see~\cite{coxeter}). Given an ordered $k$-partition $B$, we denote by $B_i$ the $i$-th set in the
$k$-tuple  $B$.

\begin{dfn}
A \emph{flag matroid} $F$ is pair $(S,\mathcal{F})$ such that $\mathcal{F}$ is
a collection of ordered $k$-partitions
of the set $S$ satisfying the following properties:
\begin{itemize}
\item[(F1)] for $1\leq i \leq k$, the set $\mathcal{B}_i=\{\cup_{1\leq j \leq i} B_j :
B\in \mathcal{F}\}$ is the set of bases of a matroid $M_i$;
\item[(F2)] for $1\leq i \leq k-1$, $M_{i}$ is a quotient of $M_{i+1}$;
\item[(F3)] if $(A_1,\ldots, A_k)$ is an ordered $k$-partition of $S$ such
that, for all $i$ with $1\leq i \leq k$, the set $A_1\cup \cdots \cup A_i$ is
a basis of the matroid $M_i$, then $(A_1,\ldots, A_k)$ is in $\mathcal{F}$.
\end{itemize}
\end{dfn}

Because of the similarity with matroids, we call the elements of $\mathcal{F}$
the \emph{flag bases} of $F$. Note that it follows from the definition
that there exist integers $r_1,\ldots, r_k$ adding up to $|S|$
such that  all ordered partitions
in $\mathcal{F}$ are in fact $(r_1,\ldots, r_k)$-partitions.
The $k$-tuple  $(r_1,\ldots, r_k)$ will be called the \emph{flag rank} of $F$.
The matroids
$M_1,\ldots M_k$ above are called the \emph{constitutents} of the flag matroid
$F$.
Notice that $M_k$ is the free matroid on $S$  and that $M_i$ has rank
$r_1+\cdots+r_i$.

A trivial example of a flag matroid is the \emph{uniform} flag matroid, having
as flag bases all possible $(r_1,\ldots,r_k)$-partitions of a set $S$.
Other examples come from chains of subspaces of a vector space, giving rise
to \emph{representable} flag matroids.
Also, given a matroid $M$, the \emph{underlying flag matroid} has as
constituents the matroids $M_i=T^i(M)$, the truncations of $M$ to ranks
$1$ to $r(M)$.
A flag matroid with flag rank~$(1,1,\ldots,1)$ can also be viewed as a
Gaussian greedoid~\cite{greedoids}.

Flag matroids, as is true of Coxeter matroids in general,  are usually viewed in
terms of their polytopes. For instance, the polytope of the uniform flag
matroid of flag rank $(1,1,\ldots,1)$ is the permutahedron; the polytope of
the underlying flag matroid is studied in~\cite{bgvw}.

%We now define matroids in terms of ordered partitions.

%\begin{dfn}\label{dfn:matroid}
%A matroid on a set $S$ is a non-empty collection $\mathcal{P}$ of
%ordered $2$-partitions of $S$ such
%that
%for any two elements $(S-B_1,B_1)$ and $(S-B_2,B_2)$ of
%$\mathcal{P}$ and any element $x\in B_1$, there exists $y\in B_2-B_1$ such
%that $(S-\{B_1\cup \{y\}}\cup \{x\}, B_1-x \cup \{y\})$ is in $\mathcal{P}$.
%\end{dfn}

%%%%%%%%%%%%%%%%%%%%%%%%%%%%%%%%%%%%%%%%%%%%%%%%%%%%%%%%%%%%%

\section{The Tennis Ball Problem}\label{sec:tbp}

%%%%%%%%%%%%%%%%%%%%%%%%%%%%%%%%%%%%%%%%%%%%%%%%%%%%%%%%%%%%%

The tennis ball problem is  a problem in enumeration that can be
phrased  in terms of balls-and-bins and in terms of lattice paths.
We need both approaches here. We first define the original problem
and show its solution amounts to counting bases of a certain type
of nested matroid. Then we generalize the problem and show that it
gives rise to a family of flag matroids.

\begin{dfn}\label{dfn:tbp}
Let $l_1$ and $l_2$ be  positive integers. Suppose we have
infinitely many balls
numbered $1,2,\ldots$ and two bins labelled $A$ and $B$. In the first
turn, balls $1,2,\ldots, l_1+l_2$ go into bin $A$, and then $l_2$ of those are
 moved
to bin $B$. In the second turn, balls $l_1+l_2+1,\ldots, 2(l_1+l_2)$ go into bin
$A$, and of the $2l_1+l_2$ balls there, $l_2$ are moved to bin $B$. At each turn,
the next $l_1+l_2$ balls go into bin $A$, and of the balls in $A$, $l_2$ are
moved to bin $B$. An \emph{$n$-configuration} is an ordered $2$-partition
of $[(l_1+l_2)n]$ giving a possible distribution of balls in the
bins after $n$ turns. The \emph{$(l_1,l_2)$-tennis ball problem} asks for the number of $n$-configurations.
\end{dfn}

This problem was solved in~\cite{tbp} using the following relationship
with nested matroids.

\begin{thm}\label{thm:tbpmats}
The number of $n$-configurations of the $(l_1,l_2)$-tennis ball problem
 is the number of bases
of the nested matroid $M[(N^{l_1} E^{l_2})^n]$.
\end{thm}

\begin{figure}[htb]
\begin{center}
\includegraphics[height=5cm]{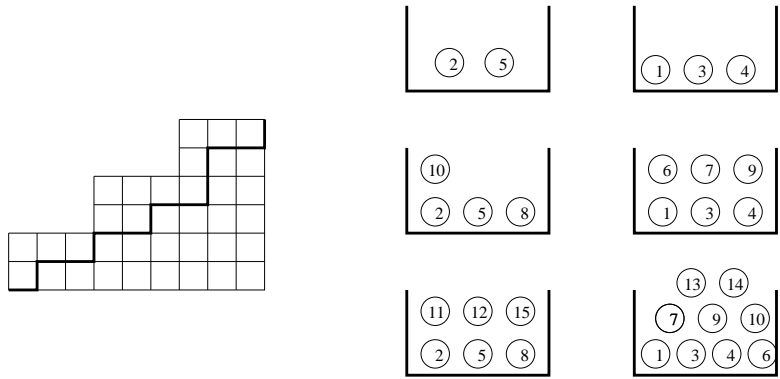} \caption{\label{fig:tbpmat} A
diagram representing the matroid $M[(N^2E^3)^3]$. The path
highlighted corresponds to the $3$-configuration shown on the
right.}
\end{center}
\end{figure}

The proof is straightforward by the bijection that sends a basis $\{n_1,\ldots,
n_r\}$ of $M[(N^{l_1}E^{l_2})^n]$ to the configuration having the balls
$\{n_1,\ldots, n_r\}$ in bin $A$ (see Figure~\ref{fig:tbpmat}). For
non-negative integers $a,b$, the matroid
$M[(N^{a} E^{b})^n]$ is called in the sequel the \emph{$n$-th $(a,b)$-tbp matroid}.

Theorem~\ref{thm:tbpmats} can be rephrased by saying that the tennis ball
problem with two bins gives the bases of a matroid. The main
result of this section is that the tennis ball problem with $k$ bins, that we
next define, gives the flag bases of a flag matroid.

\begin{dfn}\label{dfn:tbpk}
Let $(l_1,l_2,\ldots, l_k)$ be a $k$-tuple of positive integers; let
$L=l_1+l_2+\cdots +l_k$.
Suppose we have infinitely many balls
numbered $1,2,\ldots$ and $k$ bins labelled $\Gamma_1, \Gamma_2, \ldots, \Gamma_k$.
In the first turn, balls $1,2,\ldots, L$ go into bin $\Gamma_1$; of those,
 $L-l_1$
are moved
to bin $\Gamma_2$; of those, $L-l_1-l_2$ are moved to bin $\Gamma_3$, and so on
until $l_k$ balls are moved to bin $\Gamma_k$. In the second turn, balls $L+1,
\ldots, 2L$ go into bin
$\Gamma_1$, and of the $l_1+L$ balls there, $L-l_1$ are moved to bin $\Gamma_2$; of the
balls now in bin $\Gamma_2$,  $L-l_1-l_2$ are moved to bin $\Gamma_3$, and so on.
At each turn,
the next $L$ balls go into bin $\Gamma_1$, and of the balls in $\Gamma_1$, $L-l_1$ are
moved to bin $\Gamma_2$, etc... An \emph{$n$-configuration} is an ordered
$k$-partition of $[Ln]$ corresponding to a possible distribution of the
balls in the bins after $n$ turns. The \emph{$(l_1,\ldots,l_k)$-tennis ball problem} asks for the number of
$n$-configurations, that is, the number of $(nl_1,nl_2,\ldots,nl_k)$-partitions of the set $[nL]$ that
we can obtain after $n$ turns.
\end{dfn}

Note that we are only interested in
$n$-configurations, not in the movements of the balls that lead to them;
an $n$-configuration can typically be obtained by several different movements
of the balls.
Let $F_n^{(l_1,\ldots,l_k)}$ be the collection of $(nl_1,\ldots, nl_k)$-partitions of $[nL]$
that we get as $n$-configurations.   These ordered partitions are the
flag bases of a flag matroid whose constituent matroids are nested matroids.

\begin{thm}\label{thm:tbpisflag}
The set $F_n^{(l_1,\ldots,l_k)}$ is the collection of flag bases of a flag matroid on the set $[nL]$.
Moreover, for $i$ with $1\leq i \leq k$, the $i$-th constituent of the flag
matroid is the
$n$-th $(l_1+\cdots+l_i, l_{i+1}+\cdots +l_k)$-tbp matroid.
\end{thm}

\begin{proof}
We need to check that axioms (F1)--(F3) hold for
$\mathcal{F}=F_n^{(l_1,\ldots,l_k)}$. It is easy to see that
the set $\mathcal{B}_i=\{\cup_{1\leq j\leq i} B_j :
B\in \mathcal{F}\}$ is the set of bases of the $n$-th $(l_1+\cdots + l_i, l_{i+1} +
\cdots + l_k)$-tbp matroid, so (F1) holds.

To show that (F2) holds it is enough to prove
that if $a+b=a'+b'$ and $a<a'$, then the $n$-th $(a,b)$-tbp matroid $M$
is a quotient of the $n$-th $(a',b')$-tbp matroid $M'$. We show that each flat $F$ of $M$ is a flat
of $M'$. If $F$ is a cyclic flat, this follows from the characterization
of cyclic flats of nested matroids in Theorem~\ref{thm:nested}.
Otherwise, $F$ is $F'\cup I$, where $F'$ is a cyclic flat of $M$ and $I$ is
the set of isthmuses of $F$. So $F'$ is an initial segment of $[(a+b)n]$ whose
length is a multiple of $a+b$. Since $a+b=a'+b'$, by Theorem~\ref{thm:nested}
again we have that $F'$ is a
cyclic flat of $M'$. For $F'\cup I$ to be a flat of $M$, the set $I$ has to
be such that $|I\cap [t(a+b)]|< a$ for all $t$. Since $a+b=a'+b'$ and $a<a'$,
we also have that  $|I\cap [t(a'+b')]|< a'$ for all $t$, hence $F'\cup I$ is
a flat of $M'$.

To show that (F3) holds, let $(A_1,\ldots, A_k)$ be an
$(nl_1,\ldots,nl_k)$-partition of $[nL]$ such that for all $i$,
$\cup_{1\leq j\leq i} A_j$ is a basis of the $n$-th $(l_1+\cdots+
l_i, l_{i+1}+\cdots+ l_k)$-tbp matroid. We show that $(A_1,\ldots,
A_k)$ is in the collection $\mc{F}$ by showing how to get the
$n$-configuration $(A_1,\ldots, A_k)$ by suitably moving the
balls. Let $C_i$ be $\cup_{j=1}^i A_j$. To avoid wordiness, we
identify balls with the integers of their labels. We start with
all bins empty and explain how to perform $n$ turns with the
condition that at the end of each turn, the set of balls in bin
$\Gamma_i$ is a subset of $[nL]-C_{i-1}$, for all $i$ with $1\leq i \leq k$.
 Suppose we
have performed $t-1$ such turns, for $t$ with $0<t\leq n$, and let
us describe turn $t$. The assumption that
$C_{k-1},C_{k-2},\ldots,C_1$ are bases of their respective tbp
matroids gives the following facts:
\begin{itemize}
\item[$\mathbf{(k-1)}$] At least $tl_k$ of the elements of $[tL]$ are in
$[tL]-C_{k-1}$;
\item[$\mathbf{(k-2)}$] at least $t(l_k+l_{k-1})$ of the elements of $[tL]$
are in $[tL]-C_{k-2}$;\\ $\vdots$
\item[$\mathbf{(2)}$] at least $t(l_k+\cdots +l_3)$ of the elements of $[tL]$
are in $[tL]-C_2$;
\item[$\mathbf{(1)}$] at least $t(l_k+\cdots +l_2)$ of the elements of $[tL]$
are in $[tL]-C_1$.
\end{itemize}
Since after the first $t-1$ turns there are $(t-1)l_j$ elements in bin
$\Gamma_j$, we can deduce from facts $(1)$--$(k-1)$ that at this point, for all
$i$ with $1\leq i \leq k-1$,  at least $l_{i+1}+\cdots + l_k$
integers from $[tL]$ are in $[tL]-C_i$ but not in $\Gamma_{i+1}
\cup \cdots \cup \Gamma_k$.

Now move $(t-1)n+1,(t-1)n+2,\ldots,tn$ to bin $\Gamma_1$. Of all
the integers in $\Gamma_1$, choose $L-l_1$ to move to bin
$\Gamma_2$ starting by as many as possible of the ones in
$[tL]-C_{k-1}$ that are still in $\Gamma_1$, then take as many as
possible of the ones in $[tL]-C_{k-2}$, and so on until having
$L-l_1$ integers. The remarks in the previous paragraph show that
it is possible to choose integers in this way. We move them to bin
$\Gamma_2$. Since $[tL]-C_{k-1}\subset [tL]-C_{k-2}\subset \cdots
\subset [tL]-C_1$, the integers now in $\Gamma_2$ are a subset of
$[tL]-C_1$, as required. Moreover, by the way the balls are chosen, we have
that among the balls that are at this point in $\Gamma_2$, at least
$l_3+\cdots +l_k$ are in $[tL]-C_2$.

We describe generally how to move $l_i+\cdots + l_k$ balls to bin
$\Gamma_i$ from bin $\Gamma_{i-1}$ in a way such that the balls in
$\Gamma_i$ are a subset of $[tL]-C_{i-1}$, and, moreover, at least
$l_{i+1} +\cdots +l_k$ of them are in $[tL]-C_i$. From the balls in bin
$\Gamma_{i-1}$, pick $l_i+\cdots + l_{k}$ starting by as many as
possible from $[tL]-C_{k-1}$; if there are not still $l_i+\cdots +
l_{k}$, then take as many as possible from $[tL]-C_{k-2}$, and so
on, until taking as many as possible from $[tL]-C_{i-1}$. By the
same reason as above, such integers exist; the integers now in
$\Gamma_i$ are a subset of $[tL]-C_{i-1}$ and at least $l_{i+1}+\cdots
+l_k$ of them are in $[tL]-C_i$.

At the end of $n$ turns, we have $nl_i$ balls in bin $\Gamma_i$,
and these are a subset of $[nL]-C_{i-1}$, for all $i$.  This
implies that bin $\Gamma_k$ contains exactly the balls in $A_k$,
and hence bin $\Gamma_{k-1}$ contains the balls in $A_{k-1}$, and
so on. Therefore the ordered partition $(A_1, \ldots, A_k)$ is an
$n$-configuration, thus it belongs to the collection
$F^{(l_1,\ldots,l_k)}_n$ and (F3) follows.

\end{proof}

%%%%%%%%%%%%%%%%%%%%%%%%%%%%%%%%%%%%%%%%%%%%%%%%%

\section{Interpretation in terms of lattice paths}

%%%%%%%%%%%%%%%%%%%%%%%%%%%%%%%%%%%%%%%%%%%%%%%%%

The tennis ball problem with two bins has a simple interpretation in terms of
lattice paths; we associate bin $A$ with steps $N$ and bin $B$ with steps $E$,
and then each $n$-configuration corresponds to a path that does not go
above $(N^{l_1}E^{l_2})^n$. For the tennis ball problem with $k$ bins, we can
associate to each bin a direction in $\nat^k$. Then the flag bases of
$F_n^{(l_1,
\ldots, l_k)}$ are
in bijection with certain paths in $\nat^k$.
We characterize those paths combinatorially and for $k=3$ we describe them
as the set of lattice paths that do not cross a certain border.

Let $e_1,\ldots, e_k$ be the unit coordinate vectors in $\real^k$. To each
$n$-configuration
of the $(l_1,\ldots,l_k)$-tennis ball problem, we associate a path
$s_1s_2\cdots s_{nL}$ from $(0,\ldots, 0)$ to
$(nl_1,\ldots, nl_k)$ with steps defined as $s_i=e_j$ if ball $i$ is in
bin $\Gamma_j$. We call this path an \emph{$n$-configuration path}.
Hence, an $n$-configuration path can be seen as a sequence of elements from $\{e_1,
\ldots, e_k\}$. It is easy to characterize which such sequences give
configuration paths.

\begin{lemma}\label{lem:paths}
A path from $(0,\ldots, 0)$ to
$(nl_1,\ldots, nl_k)$ is an $n$-configuration path for the
$(l_1,\ldots,l_k)$-tennis
ball problem if and only if, for all $t$ with $1\leq t \leq n$
and all $i$ with $1\leq i \leq k-1$, among the first $tL$ steps there are
at most $t(l_1+\cdots+l_i)$ whose type belongs to $\{e_1,\ldots,e_i\}$.
\end{lemma}

\begin{proof}
After $t$ turns, for $1\leq t \leq n$, there are exactly $t(l_1+\cdots +l_i)$
balls of  $[tL]$ in the first $i$ bins. Since balls can only move to bins
with a higher index, at the end of $n$ turns there are at most $t(l_1+\cdots
+l_i)$ balls of $[tL]$ in $\Gamma_1 \cup \cdots \cup \Gamma_i$.
Hence among the first $tL$
steps of a configuration path there are at most $t(l_1+\cdots+l_i)$ steps
whose type is in $\{e_1,\ldots, e_i\}$.

For the converse, assume we have a path $\pi$ that satisfies
the condition. Let $A_i$ be the set of integers $s$ such that step $s$ in $\pi$ is of type $e_i$. Consider the $(nl_1,\ldots, nl_k)$-partition
$(A_1,\ldots, A_k)$ of $[nL]$ obtained in this way. The condition on the path
implies that the set $\cup_{j=1}^i A_j$ is a basis of the $(l_1+\ldots +l_i,
l_{i+1}+\ldots + l_k)$-tbp matroid, for all $i$ with $1\leq i \leq k-1$.
Therefore $(A_1,\ldots, A_k)$ is a flag basis of the flag matroid $F_{n}^{(l_1,
\ldots, l_k)}$ and hence $\pi$ is an $n$-configuration path, as required.
\end{proof}

The following is an immediate corollary.

\begin{cor}\label{cor:movepaths}
Given an $n$-configuration path $\pi$, the path obtained by
switching a pair of steps $s_i=e_l$ and $s_j=e_m$ is a
configuration path if $i<j$ and $l\leq m$. Moreover, let $\pi'$ be
an initial segment of $\pi$ with $t'_j$ steps of type $e_j$
for all $1\leq j\leq k$. Let
$n'$ be the minimum integer such that $t'_j\leq n'l_j$ for all $j$.
For $n''\geq 0$, consider the path obtained from $\pi'$ followed by
$(n'+n'')l_k-t'_k$ steps $e_k$, then $(n'+n'')l_{k-1}-t'_{k-1}$ steps
$e_{k-1}$, and so on, until finishing with $(n'+n'')l_1-t'_1$ steps
$e_1$. Then this path is an $(n'+n'')$-configuration path.
\end{cor}

The \emph{$n$-diagram} for the $(l_1,\ldots, l_k)$-tennis ball problem
  is the set of points in $\nat^k$ that are contained in some
$n$-configuration path. Our goal is to study what these diagrams look like
for $k=3$. The
following is another corollary of Lemma~\ref{lem:paths}.

\begin{cor}\label{cor:push}
If $(x,y,z)$ is in the $n$-diagram, then $(x',y,z')$ is in the $n$-diagram
for all
$z\leq z' \leq nl_3$ and all $x'\leq x$.
\end{cor}

A direct consequence of this corollary is that to describe the $n$-diagram it is
enough to give, for each $(x,y)$, the minimum value of $z$ such that
$(x,y,z)$ is a point of the $n$-diagram; this minimum is denoted $m_n(x,y)$.
(Trivially the maximum value of $z$ is
$nl_3$.)
If no such $z$ exists, because no point of the form $(x,y,\ast)$ is in the
$n$-diagram, we set $m_n(x,y)=\ast$.
So the $n$-diagram is described
by an $(nl_1+1) \times (nl_2+1)$ matrix $\mc{M}_n$ with entries in the set
$\{0,1,\ldots, nl_3\}\cup \{\ast\}$ and such that in row~$x$ and column~$y$
we have $m_n(x-1,y-1)$. If $n=1$, then
trivially $\mc{M}_1$ is the zero matrix and the corresponding $1$-diagram is represented in
Figure~\ref{fig:n1_243}. In all figures below, the direction of the third
coordinate has been reversed for a better view of the picture, and, as pointed
out above, all points under a point that is shown are in the diagram as well.

\begin{figure}[htb]
\begin{center}
\includegraphics[height=5cm]{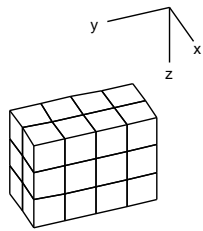} \caption{\label{fig:n1_243}
The $1$-diagram for the $(2,4,3)$-tennis ball problem.}
\end{center}
\end{figure}

\begin{figure}[htb]
\begin{center}
\includegraphics[height=8cm]{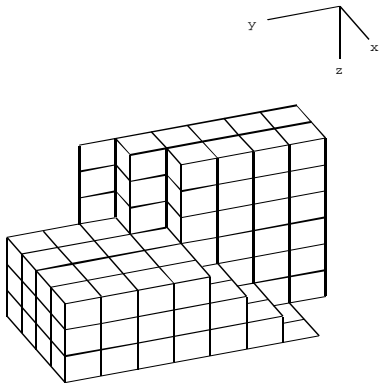} \caption{\label{fig:n2_243}The $2$-diagram
for the $(2,4,3)$-tennis ball problem.}
\end{center}
\end{figure}

The $2$-diagram for the $(2,4,3)$-tennis ball problem is shown in
Figure~\ref{fig:n2_243}. We first give the matrices $\mc{M}_n$ and then prove
that they give
the right diagrams.

The matrix $\mc{M}_2$ is
made up of four blocks,
$$\mc{M}_2=\left(\begin{array}{c|c} A & B\\ \hline  C&D \\\end{array}\right),$$
where $A$ is the $(l_1+1)\times (l_2+1)$ matrix all whose entries
are zero (hence, it is $\mc{M}_1$), $D$ is an $l_1\times l_2$
matrix all whose entries are $l_3$, $B$ is an $(l_1+1)\times l_2$
matrix with
$$b_{i,j}=\left\{\begin{array}{ll}0& \mathrm{if }\ i\leq l_1+1-j,\\ l_3 &
\mathrm{otherwise},\end{array} \right. $$
and $C$ is the $l_1\times (l_2+1)$ matrix with $\min\{l_2,l_3\}+1$ non-$\ast$
columns, with the elements in the last column being $l_3$ and each other
non-$\ast$ column being obtained by adding $+1$ to the next, that is,

$$ \left(\begin{array}{ccccccccc} \ast & \cdots & \ast &2l_3 & 2l_3-1& 2l_3-2& \cdots&
l_3+1&l_3\\ \vdots & \ddots &\vdots &\vdots &\vdots &\vdots &\ddots & \vdots& \vdots\\
\ast & \cdots & \ast &2l_3 & 2l_3-1& 2l_3-2& \cdots&
l_3+1&l_3\\\end{array}\right).$$

\begin{figure}[ht]
\begin{center}
\includegraphics[width=0.8\textwidth]{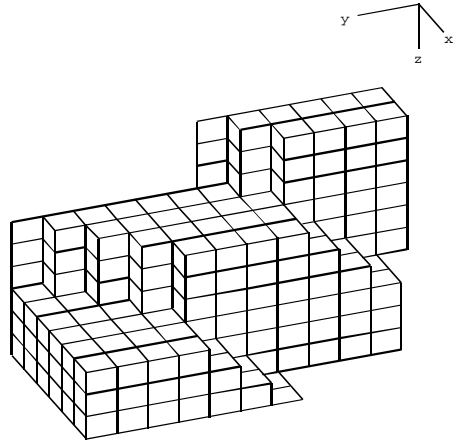}
\caption{\label{fig:n3_243}The $3$-diagram for the $(2,4,3)$-tennis ball
problem.}
\end{center}
\end{figure}

\begin{figure}[ht]
\begin{center}
\includegraphics[width=5cm]{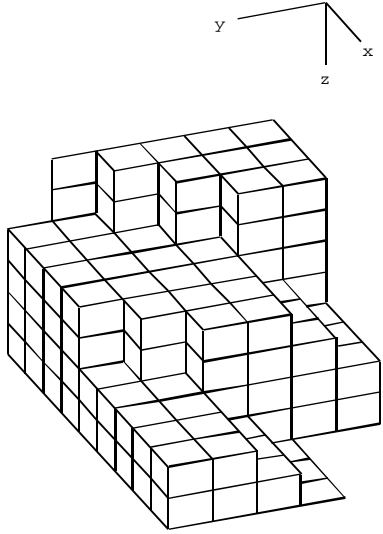}
\caption{\label{fig:n3_322}The $3$-diagram for the $(3,2,2)$-tennis ball
problem.}
\end{center}
\end{figure}

We now define recursively the $(nl_1+1)\times (nl_2+1)$ matrix
$\mc{M}_n$ that gives
 the $n$-diagram. The matrix $\mc{M}_n$ also decomposes into $4$
 blocks
 $$\mc{M}_n=\left(\begin{array}{c|c} A_n & B_n\\ \hline  C_n&D_n \\\end{array}\right),$$
where $A_n$ has $(n-1)l_1+1$ rows and $(n-1)l_2+1$ columns. The
matrix $A_n$ is $\mc{M}_{n-1}$. The entry in row $i$ and column
$j$ of the matrix $B_n$, for $1\leq i \leq (n-1)l_1+1$ and $1\leq
j \leq l_2$, is given by $(n-s)l_3$, where $s$ is the only integer
for which $(n-s)(l_1+l_2) < i-1+j+(n-1)l_2\leq (n-s+1)(l_1+l_2)$.
%$$l_3\left(n-1-\lceil \frac{(n-1)l_1-i-j+2}{l_1+l_2} \rceil\right) $$
%$(n-k)l_3$ for the only integer $k$ such that
%$(n-k)(l_1+l_2)<(i+1)+(j+1)\leq (n-k+1)(l_1+l_2)$.
Roughly
speaking, $B_n$ consists of diagonal stripes of width $l_1+l_2$, see the examples below.
 The matrix $C_n$ has
$l_1$ rows and $(n-1)l_2+1$ columns; all entries in the last
column are $(n-1)l_3$ and each column is obtained by adding
one to the next, until we reach $nl_3$; hence it  is given by
$$ \left(\begin{array}{ccccccccc} \ast & \cdots & \ast &nl_3 & nl_3-1& nl_3-2& \cdots&
(n-1)l_3+1&(n-1)l_3\\ \vdots & \ddots &\vdots &\vdots &\vdots &\vdots &\ddots & \vdots& \vdots\\
\ast & \cdots & \ast &nl_3 & nl_3-1& nl_3-2& \cdots&
(n-1)l_3+1&(n-1)l_3\\\end{array}\right).$$ Finally, $D_n$ is the
$l_1\times l_2$ matrix all whose entries are $(n-1)l_3$. The
matrix $\mc{M}_3$ is shown below for the $(2,4,3)$- and the $(3,2,2)$-tennis
ball problems, and the corresponding $3$-diagrams are shown in Figures~\ref{fig:n3_243} and~\ref{fig:n3_322}.

$$
\left(\begin{array}{c} \begin{array}{ccccc|cccc|cccc}
0&0&0&0&0&0&0&3&3&3&3&3&3  \\
0&0&0&0&0&0&3&3&3&3&3&3&6          \\
0&0&0&0&0&3&3&3&3&3&3&6&6          \\
\cline{1-9}
\ast&6&5&4&3&3&3&3&3&3&6&6&6         \\
\ast&6&5&4&3&3&3&3&3&6&6&6&6       \\
\hline\end{array}\\ \begin{array}{ccccccccc|cccc}
\ast&\ast&\ast&\ast&\ast&9&8&7&6&6&6&6&6 \\
\ast&\ast&\ast&\ast&\ast&9&8&7&6&6&6&6&6 \\
\end{array}\end{array}\right) \quad \left(\begin{array}{c}\begin{array}{ccc|cc|cc}
0&0&0&0&0&0&2\\
0&0&0&0&0&2&2\\
0&0&0&0&2&2&2\\
0&0&0&2&2&2&2\\
\cline{1-5}
4&3&2&2&2&2&2\\
4&3&2&2&2&2&4\\
4&3&2&2&2&4&4\\
\hline\end{array}\\
\begin{array}{ccccc|cc}
\ast&\ast&6&5&4&4&4\\
\ast&\ast&6&5&4&4&4\\
\ast&\ast&6&5&4&4&4\\
\end{array}\end{array}  \right) $$

We now show that the matrix $\mc{M}_n$ gives the $n$-diagram and, moreover, that all paths
contained in the $n$-diagram are $n$-configuration paths.

\begin{thm}\label{thm:diagram}
The $n$-diagram for the $(l_1,l_2,l_3)$-tennis ball problem is given by the
matrix $\mc{M}_n$. Furthermore, the $n$-configuration paths are exactly those
contained in the $n$-diagram.
\end{thm}

\begin{proof}
The proof is by induction on $n$. As seen above the case $n=1$ is trivial.
Assume $\mc{M}_{n-1}$ is the matrix of $(n-1)$-diagram for the
$(l_1,l_2,l_3)$-tennis ball problem. Let $\mc{N}$ be the matrix
of the $n$-diagram; we prove that $\mc{N}=\mc{M}_n$. The matrix
$\mc{N}$ has dimensions $(nl_1+1)\times (nl_2+1)$.

Recall that the entry in row $x+1$ and column $y+1$ of $\mc{N}$ is
$m_n(x,y)$.  Given $x$ and $y$ with $0\leq x\leq nl_1$ and $0\leq
y \leq nl_2$, let $z$ be $m_n(x,y)$. We show that $z$ is the entry
in row $x+1$ and column $y+1$ of the matrix $\mc{M}_n$. The proof has
 three cases.

\medskip

\textbf{Case 1.} $x\leq (n-1)l_1$ and $y\leq (n-1)l_2$

\smallskip

In this case we need to show that $z=m_{n-1}(x,y)$.
Let $\pi$ be an $(n-1)$-configuration path that contains the point
$(x,y,m_{n-1}(x,y))$. By Corollary~\ref{cor:movepaths}, $\pi$ can be
extended to an $n$-configuration path, hence $z \leq
m_{n-1}(x,y)$. To show equality, assume $z< m_{n-1}(x,y)\leq
(n-1)l_3$
and let $\rho$ be an $n$-configuration path that contains the
point $(x,y,z)$. Let $\rho'$ be the initial segment corresponding
to the first $x+y+z$~steps. By Corollary~\ref{cor:movepaths} again, $\rho'$
can be extended to an $(n-1)$-configuration path
contradicting the induction hypotheses.

\medskip

\textbf{Case 2.} $x>(n-1)l_1$

\smallskip

We want to show in this case that $z$ is given by the entries of
the matrices $C_n$ and $D_n$, depending on the value of $y$.
Observe first that Lemma~\ref{lem:paths} implies that $x+y+z>
(n-1)L$ and hence also that $y+z \geq (n-1)(l_2+l_3)$.

We discuss now the subcase $y\geq (n-1)l_2$,  showing that
$z=(n-1)l_3$. Let $\pi$ be any $(n-1)$-configuration path. Extend
$\pi$ by adding $y-(n-1)l_2$ steps $e_2$ followed by $x-(n-1)l_1$
steps $e_1$ and then $l_3$ steps $e_3$, and add the remaining steps in any
way. This path clearly satisfies the
condition on Lemma~\ref{lem:paths}, hence $z\leq (n-1)l_3$. To
complete the proof of the claim, suppose that $z<(n-1)l_3$. Then
by Lemma~\ref{lem:paths} $x+y+z<(n-1)L$ . This contradicts the
first conclusion of the
previous paragraph.

The other subcase left is $y<(n-1)l_2$. We show that in this case
$z=(n-1)(l_2+l_3)-y$. Since we already know that $z$ is at least
$(n-1)(l_2+l_3)-y$ it is enough to show that there is an
$n$-configuration path containing the point
$(x,y,(n-1)(l_2+l_3)-y)$. Consider the path
$$\sigma =
(e_3)^{(n-1)(l_2+l_3)-y}  (e_2)^y
 (e_1)^x (e_3)^{nl_3-(n-1)(l_2+l_3)+y}
 (e_2)^{nl_2-y}  (e_1)^{nl_1-x}.$$
 It is easy to check that this path satisfies the condition of
Lemma~\ref{lem:paths} and hence it is an $n$-configuration path
containing the point $(x,y,(n-1)(l_2+l_3)-y)$, as required.

\medskip

\textbf{Case 3.} $x<(n-1)l_1$ and $y>(n-1)l_2$

\smallskip

In this case the value of $z$ has to be the one given by the matrix $B_n$.
We have to show that if  $s$ is such that
 $(n-s)(l_1+l_2)\leq x+y <(n-s+1)(l_1+l_2)$ then
$z=(n-s)l_3$.

Since from Case~2 we have that the point $(nl_1,y,(n-1)l_3)$ is in the
$n$-diagram, by
Corollary~\ref{cor:push} it follows that the point $(x,y,(n-1)l_3)$
is in the $n$-diagram as well, hence $z\leq (n-1)l_3$.
If $x+y> (n-1)(l_1+l_2)$, then by Lemma~\ref{lem:paths}
$x+y+z>(n-1)L$ and hence $z\geq (n-1)l_3$, so in this case $z=(n-1)l_3$.

Now suppose $(n-s)(l_1+l_2)< x+y \leq (n-s+1)(l_1+l_2)$. As in the
previous paragraph, Lemma~\ref{lem:paths} implies that $z\geq
(n-s)l_3$. To show that we have equality, we give an
$n$-configuration path containing the point $(x,y,(n-s)l_3)$.
Consider the path
$$\sigma= (e_3)^{(n-s)l_3} (e_2)^y (e_1)^x (e_3)^{sl_3} (e_2)^{nl_2-y} (e_1)^{nl_1-x}. $$
By using that $y>(n-1)l_2$ and that $x+y\leq (n-s+1)(l_1+l_2)$ it
is easy to show that $\sigma$ satifies the condition of
Lemma~\ref{lem:paths} and hence it is an $n$-configuration path.

\medskip

To finish the proof we have to show that any path contained in the $n$-diagram
 is an $n$-configuration path. Let $\pi$ be such a path; we check that $\pi$
satisfies the condition in Lemma~\ref{lem:paths}. Let $(X,Y,Z)$ be a point
in the path with $X+Y+Z=tL$ for some $t$ with $1\leq t\leq n-1$; our goal is
to show that $X\leq tl_1$ and $X+Y \leq t(l_1+l_2)$. Consider the point
$p=(X,Y,m_n(X,Y))$. The proof above shows that there is an $n$-configuration
path $\pi'$ that goes through $p$; since $Z\geq m_n(X,Y)$, we can apply
Corollary~\ref{cor:movepaths} to obtain from $\pi'$ an $n$-configuration
path containing the point $(X,Y,Z)$. Since all $n$-configuration paths
satisfy the condition in Lemma~\ref{lem:paths}, we have that $X\leq tl_1$ and
$X+Y\leq t(l_1+l_2)$.
\end{proof}

\section{Concluding remarks}

\begin{figure}[ht]
\begin{center}
\includegraphics[width=3cm]{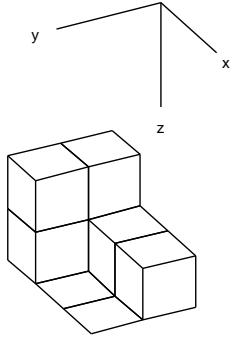}
\caption{\label{fig:noflag}The set of paths contained in the diagram does not give rise to a
flag matroid.}
\end{center}
\end{figure}

The results of the previous section show that some sets of lattice paths in
$3$ dimensions can be interpreted in terms of flag matroids; hence, the flag
matroids one obtains from the tennis ball problem naturally generalize lattice
path matroids.
This might lead to the suspicion that any set of paths in $\nat^k$ with a
``reasonable'' border also gives rise to flag matroids. Unfortunately,
it is very easy to
produce counterexamples to this. For instance, consider the diagram in
Figure~\ref{fig:noflag}. If the paths contained in that diagram were in correspondence with
the flag bases of a flag matroid $F$, we would have that $B_1=\{2,6\}$
and $B_2=\{4,5\}$ are cobases of the second constituent of $F$.
Hence  it should be possible to replace $2$ in $B_1$
by either $4$ or $5$. But no path contained in the diagram has $\{4,6\}$ or
$\{5,6\}$ as its set of steps in the direction $e_3$.

A question that remains open is to solve the $(l_1,\ldots, l_k)$-tennis
ball problem, or even the $(1,1,1)$-tennis ball problem.  The approaches
used previously to solve the case $k=2$ do not seem to
generalize easily. In particular, the strategy from~\cite{tbp} would suggest
the use of a Tutte polynomial-like invariant to count
flag bases. There are some
generalizations of the Tutte polynomial to pairs of matroids and to chains of
matroids related by
strong maps (\cite{lv,wk}), but unfortunately they do not seem to include the
number of flag
bases as an specialization. Following the Tutte polynomial approach for flag matroids would require
defining first the suitable generalization.

We finish with some easy bounds. A trivial upper
 bound for
 the number of $n$-configurations is given by the
total number of $(nl_1,\ldots,nl_k)$-partitions of $[nL]$,
which is the multinomial
coefficient $$\binom{nL}{nl_1,\ldots,nl_k}=\frac{(nL)!}{(nl_1)!\cdots (nl_k)!}.$$
%\sim
% C\left(\frac{(l_1+\cdots +l_k)^{l_1+\cdots +l_k}\right)^n
%n^{-\frac{k-1}{2}}, $$
%for some constant $C$ depending on $l_1,\ldots,l_k,$ and $k$
%and by using Stirling's approximation of the factorial.
The following connection with
 Young tableaux gives a lower bound on the number of $n$-configurations of
the $(l_1,\ldots,l_k)$-tennis ball problem when $l_k\geq l_{k-1}\geq \cdots
\geq l_1$. Consider sequences of length $n(l_1+\cdots l_k)$
over the alphabet $\{e_1,\ldots,e_k\}$ containing $nl_i$ copies of $e_i$ and
 such that in any
initial subsequence the number of symbols $e_i$ is greater or equal than the
number of symbols $e_{i-1}$, for all $i$ with $2\leq i \leq k$. Since all
these sequences trivially satisfy the condition in Lemma~\ref{lem:paths}, they
give $n$-configuration paths. The number of such sequences equals the number
of standard Young tableaux of shape $(nl_k,\ldots,nl_1)$, and this is given by
the hook-length formula (see~\cite[Chapter 7]{ec2}). In the case
$l_1=\cdots=l_k=l$, this is $$\frac{(nlk)!}{(nl)!^k \prod_{i=1}^{k-1}(\frac{nl}{i}+1)^{k-i}}.$$ The general case gets more involved and we omit it since not
much insight is gained.

For the $(1,1,1)$-tennis ball problem, the first order approximation
of the lower and upper bounds are $C 27^n n^{-7/2}$ and $C' 27^n n^{-1/2}$,
respectively, for some constants $C$ and $C'$. Computer evidence seems to
suggest that the right number lies closer to the lower bound, and that the
 exponent in the term on $n$ is $-3$ (\cite{dominic}).

A general lower bound for the $(l_1,\ldots, l_k)$-tennis ball problem can
be obtained as follows. Let $t(a,b,n)$ be the number of $n$-configurations
of the $(a,b)$-tennis ball problem. Then the number of $n$-configurations
of the $(l_1,\ldots, l_k)$-tennis ball problem is at least
$$t(l_1,l_2+\cdots +l_k,n)t(l_2,l_3+\cdots +l_k,n)\cdots t(l_{k-1},l_k,n), $$
since we can think of the $(l_1,\ldots, l_k)$-tennis ball problem as $n$
turns of the
$(l_1,l_2+\cdots +l_k)$-tennis ball problem, followed by $n$ turns of the
$(l_2,l_3+\cdots +l_k)$-tennis ball problem on the result of the first, and so on.
The bound is strict since each $t(a,b,n)$ counts the number of
$n$-configurations of the $(a,b)$-tennis ball problem, but each of these can
usually be reached by several movements of the balls and that is relevant
for the version with $k$ bins.

%whose first term in the asymptotic
%expansion is $$C' k^{lkn} n^{-\frac{k-1}{2}-\binom{k}{2}},$$ for some
%constant $C'$ depending on $k$ and $l$
%$$\frac{(nL)!}{\Prod_{i,j} (n\lambd_i +\lambda'_j-i-j+1} $$

\section*{Acknowledgements}

The author thanks Dominic Welsh for asking a question
that suggested this work, and Joseph Bonin and Sergi Elizalde for helpful
discussions.

\end{document}